\emergencystretch0.4mm
\input AHTOH-E.STY
\def\CODIM{\={\rm codim\,}}
\def\C{{\cal C}}
\def\V{{\cal V}}
\def\N{{\cal N}}

\def\P{{\cal P}}
\def\F{{\cal F}}
\def\Ends{{\rm Ends\,}}

%%%%%%%%%%%%%%%%%%%%%%%%%%%%%%%
\UDC{\let\,\relax\kern-3pt\
512.544.44\,+\,512.543.22%
\,+\,512.552.16\,+\,512.552.4}
\MSC{20e36, 20e10%
, 17a30, 17a36}
\title{Automorphism invariance and identities}
\author{%
{Evgenii I. Khukhro$^\sharp$}
\quad
{Anton A. Klyachko$^\flat$}
\quad
{Natalia Yu. Makarenko$^{\sharp\natural}$}
\quad
{Yulia B. Melnikova$^\flat$}
}
\address{
\\
$^\sharp$Sobolev Institute of Mathematics, Novosibirsk 630090,
prospekt akademika Koptyuga 4.
\\
$^\flat$Faculty of Mechanics and Mathematics, Moscow State University,
Moscow 119991, Leninskie gory, MSU.
\\
$^\natural$Universit\'e de Haute Alsace (UHA)
\\
{khukhro@yahoo.co.uk}
\quad
{klyachko@mech.math.msu.su}
\quad
{makarenk@math.nsc.ru}
\quad
{yuliamel@mail.ru}
}
\footnote{}{\hskip-0.5cm The work of the second author was
supported by \RFBR 08-01-00573.}
\footnote{}{\hskip-0.5cm The work of the third author was
supported by the Council on Grants of the President of the Russian
Federation, program of support of the leading scientific schools, grant
NSh-344.2008.1.}
\abstract{%
If an outer (multilinear) commutator
identity holds in a large subgroup of
a group, then it holds also in a
large
characteristic subgroup. Similar assertions are valid for
algebras and their ideals or subspaces. Varying the
meaning of the
word ``large", we obtain many interesting facts.
These results cannot be extended to arbitrary
(non-multilinear) identities.
As an application, we give a
sharp estimate
for the `virtual derived length' %!!
of
(virtually solvable)-by-(virtually solvable) groups.
}

%%%%%%%%%%%%%%%%%%%%%%%%%%%%%%%%%%%%%%%%%%%%%%%%%%%%%%%%%%%%%%%%%%%%%%%%%%
%%%%%%%%%%%%%%%%%%%%%%%%%%%%%%%%%%%%%%%%%%%%%%%%%%%%%%%%%%%%%%%%%%%%%%%%%%
\s 0.
Introduction

In [KhM07a], it was shown that each group virtually satisfying an
outer commutator identity contains a finite-index characteristic
subgroup satisfying this identity. A similar result was obtained
in [KhM08] for ideals in arbitrary algebras (with codimension
taking the role of index). Another similar result can be found in
[KhM07b]: if a finite $p$-group $G$ contains a normal class $t$ %!!
nilpotent subgroup $N$, then it also contains a characteristic
class $t$ nilpotent subgroup $H$ whose co-rank is bounded by a
function of the co-rank of $N$ and the number $t$. The paper
[KlM09] contains a new much shorter proof of the theorem about
index and a better estimate for the index of the characteristic
subgroup.

 In this note, we extend and generalize the new proof in
 [KlM09] to a wide class of algebraic systems. In this
 generalization the property of having finite index is replaced by
 a certain abstract property of ``smallness'', and index is
 replaced by a certain abstract ``codimension''. This general
 theorem includes all the aforementioned results as special cases.
 Moreover, there are many other new applications, which include,
 for example, the case of finite $p$-groups with subgroups of
 bounded co-rank satisfying arbitrary outer commutator
 identity (rather than just nilpotency identity as in
 [KhM07b]).

 On the other hand, we show that the original proofs in
[KhM07a], [KhM07b], and [KhM08] yield a bit more. Namely, each
group contains only finitely many %!!
finite-index subgroups which are maximal (by inclusion) among all
normal %!!
subgroups satisfying a given outer commutator identity. In
particular, this implies that each finite-index subgroup which is
maximal among all normal %!!
subgroups satisfying a given outer
commutator identity contains a finite-index characteristic
subgroup. Similar stronger results are valid for algebras over
fields.

To complete the picture, we mention earlier known results on this
subject. Let $G$ be a group and let $N$ be its finite-index
subgroup. Then
\- $N$ contains a normal in $G$ subgroup of finite
index (dividing $|G{\,:\,}N|!$);
\- if $G$ is finitely generated,
then $N$ contains a fully %!!
invariant (and even verbal) in $G$ subgroup of finite index;
\- if
$N$ is abelian, then $G$ contains a characteristic abelian
subgroup of finite index.

\noindent
These facts are well known and can be found in textbooks
on group theory (see, e.g., [KaM82]). Note also that in [BeK03] it
was proved that the existence of a solvable finite-index subgroup
of derived length %!!
$t$ implies the existence of a characteristic
solvable finite-index subgroup of derived length $\le t^2$.
In [BrNa04], it was shown that 
any virtually nilpotent group contains a characteristic nilpotent
finite-index subgroup.

\smallskip

As an application of the obtained results, in Section 5 we obtain
a sharp estimate for the `virtual derived length' of extensions of
virtually solvable groups by virtually solvable groups. This
answers a question of J.~Button.

In Section~6, we consider the periodicity law $x^p=1$. We show
that the theorem about finite-index characteristic subgroup
does not hold for this identity (if $p$ is a large prime).

%%%%%%%%%%%%%%%%%%%%%%%%%%%%%%%%%%%%%%%%%%%%%%%%%%%%%%%%%%%%%%%%%%%%%%%
\s 1.
The results

\Th 1 \rm[KhM07a], [KlM09]. If a group $G$ contains a normal
finite-index subgroup $N$ satisfying an outer commutator identity
$w(x_1,\dots,x_t)=1$, then $G$ contains a characteristic and even
invariant under %!!
all surjective endomorphisms subgroup $H$ satisfying the same
identity and such that
$
\log_2|G{\,:\,}H|\le f^{t-1}(\log_2{|G{\,:\,}N|}).
$

Henceforth, $f^k(x)$ means the $k$-th iteration of the function
$f(x)=x(x+1)$. An \emph{outer \({\rm or} multilinear\) commutator
identity} is an identity of the form
$[\dots[x_1,\dots,x_t]\dots]=1$ with some meaningful arrangement
of brackets, where all letters $x_1,\dots,x_t$ are different.
Examples of such identities are solvability, nilpotency,
centre-by-metabelianity, etc. A formal definition looks as
follows. Let $F(x_1,x_2,\dots)$ be a free group of countable rank.
An \emph{outer commutator of weight 1} is just a letter $x_i$. An
\emph{outer commutator of weight $t>1$} is a word of the form
$w(x_1,\dots,x_t)=[u(x_1,\dots,x_r),v(x_{r+1},\dots,x_t)]$, where
$u$ and $v$ are outer commutators of weights $r$ and $t-r$,
respectively. An \emph{outer commutator identity} is an identity
of the form $w=1$, where $w$ is an outer commutator.

\Remark 1. The condition that the subgroup $N$
be
normal is not
essential. It is well known that any finite-index subgroup $N$
contains a normal finite-index subgroup $\~N$ such that
$|G{\,:\,}\~N|$ does not exceed (and even divides) $|G{\,:\,}N|!$
(see, e.g., [KaM82]). Therefore, Theorem 1 remains valid for
non-normal subgroups $N$, but with worse estimate
$\log_2|G{\,:\,}H|\le f^{t-1}(\log_2{|G{\,:\,}N|!})$.

\Remark 2.
Theorem $1'$ (see below) implies that, under the conditions of Theorem 1,
a characteristic (and even invariant under all surjective
endomorphisms) finite-index subgroup can be found inside any
finite-index subgroup which is maximal (by inclusion) among all normal
subgroups satisfying the outer commutator identity.

\Remark 3.
Theorem 4 (see below), which is a generalization of Theorem 1,
implies that the group~$G/H$ lies in the variety generated by
$G/N$ (and even in the formation generated by this group).

\Th 2 {\rm(cf. [KhM08])}. Let $G$ be an algebra (possibly,
non-associative) over a field. If $G$ contains a
finite-codimensional subspace $N$ satisfying a multilinear
identity $w(x_1,\dots,x_t)=0$, then $G$ contains a subspace $H$
satisfying the same identity, invariant under all surjective
endomorphisms, and such that $\codim H \le f^{t-1}(\codim N)$.
This subspace $H$ is left, right, or two-sided ideal if the
subspace $N$ is left, right, or two-sided ideal, respectively.

\Remark. Theorem $2'$ (see below) implies that, under the
conditions of Theorem 2, a finite-codimensional subspace invariant
under all surjective endomorphisms can be found inside any
finite-codimensional subspace which is maximal (by inclusion)
among all subspaces satisfying the multilinear identity. A similar
fact is valid for ideals (left, right, and two-sided).

\Th 3 {\rm(cf. [KhM07b])}.
If a finite $p$-group $G$ contains a normal subgroup $N$
satisfying an outer commutator identity
$w(x_1,\dots,x_t)=1$, then $G$ contains a
characteristic subgroup~$H$
satisfying the same identity and such that
$
\rank G/H \le
f^{t-1}(\rank G/N).
$

Here, $\rank G$ is the minimal positive integer $n$ such that
any finitely generated subgroup of $G$ is generated by at most $n$
elements.

\Th 4. If a group $G$ contains a normal subgroup $N$ satisfying an
outer commutator identity \hbox{$w(x_1,\dots,x_t)=1$} and $G/N$ has a
smallness property $\P$, then $G$ contains a characteristic and
even invariant under all surjective endomorphisms subgroup $H$,
satisfying the same identity and such that $G/H$ has the property
$\P$.

\noindent A \emph{smallness property} in Theorem 4 is any abstract
group property $\P$  %!!
satisfying the following conditions:
\item {1)}
a quotient group of a group with property $\P$ also has this property;
\item {2)}
a subdirect product of two groups with property $\P$ also has
this property;
\item {3)}
each group with property $\P$ satisfies the maximality (ACC) condition
for normal subgroups.

\enditem
Examples of such properties are the maximality condition,
the maximality condition for normal subgroups, polycyclicity, finiteness,
etc.

Theorem 1 was proved in [KhM07a], but with a worse estimate for the
index. A simpler proof and the estimate presented above was obtained in
[KlM09]. Theorem 2 was proved in [KhM08] with a worse estimate for
the codimension. An important particular case of Theorem 3 corresponding
to the nilpotency identity was proved in [KhM07b]. Theorem~4 is new. All
these theorems turn out to be special cases of a general fact
concerning multi-operator groups.

The methods developed in [KhM07a], [KhM07b], and [KhM08] allows us to
prove assertions stronger than Theorems 1 and 2, but with
worse estimates.

\Th 1$'$.
Let $w(x_1,\dots,x_t)$ be an outer commutator. Then,
in any group, the number of finite-index subgroups which are maximal
(by inclusion) among all normal subgroups satisfying
the identity $w(x_1,\dots,x_t)=1$ is finite.
Moreover, the number of such subgroups of index $\le n$ does not exceed
$$
2^{F^{t-1}(n)},
\quad\hbox{
where $F^k(x)$ is the $k$-th iteration of the function $F(x)=xn^{2^x}$.}
$$

\Remark. This theorem consists of two independent assertions. On
the one hand, the number of subgroups of index $\le n$ which are
maximal among all normal subgroups with the identity
$w(x_1,\dots,x_t)=1$ is bounded by an explicit function of $n$ and
$t$. This function grows very fast, but on the other hand, the
total number of finite-index subgroups which are
maximal among all normal
%!!!!!!!???? russkii!!
%??
subgroups with given
identity is finite. The following theorem shows that a similar
statement is valid for subspaces (or ideals) in algebras.

\Th 2$'$.
Let $w(x_1,\dots,x_t)$ be a multilinear element of the free
(non-associative) algebra over a field $F$. Then,
in any algebra over $F$, the intersection of all
finite-codimensional ideals which are
maximal (by inclusion) among all ideals satisfying the identity
$w(x_1,\dots,x_t)=0$ has finite codimension.  Moreover, the intersection
of such ideals of codimension $\le n$ has codimension not larger
than some number depending only on $n$ and $t$. Here, the word
``ideal" means left, right, two-sided ideal, or simply subspace (zero-side
ideal).

Theorem 1 makes it possible to give the most exact answer to a question of
J.~O.~Button ([But08], Problem~3) on extensions of virtually solvable
groups by virtually solvable groups.

\Th 5. Any extension of a virtually solvable of derived length $s$
group by a virtually solvable of derived length~$t$ group is
virtually solvable of derived length $\le t+s+1$.

In section 5, we prove this theorem and give a simple example
showing that the obtained estimate is
best possible.

The following assertion shows that neither Theorem 1 nor other
versions of Theorem 4 can be extended to arbitrary identities.

\Th 6.
For any sufficiently large prime $p$, there exists a
group $G$ of exponent $p^2$
with a finite-index subgroup of exponent $p$, but without
characteristic finite-index subgroups of exponent $p$. Moreover,
any quotient of $G$ by a characteristic subgroup
of exponent $p$ satisfies neither the maximality condition for normal
subgroups nor the minimality condition for normal subgroups and,
therefore, has no smallness properties.

%%%%%%%%%%%%%%%%%%%%%%%%%%%%%%%%%%%%%%%%%%%%%%%%%%%%%%%%%%%%%%%%%%%%%%%%%%
\s 2.
Multi-operator groups

Recall that, according to [Kur62], an \emph{$\Omega$-group} is a group
$(G,+)$ (not necessarily commutative) with a family of operations
$\Omega$. Each operation $f\in\Omega$ is a mapping
$f:G^{n_f}\to G$ from a finite Cartesian power of $G$ to $G$ such that
$f(0,\dots,0)=0$.

\noindent{\bf Examples.}

\nobreak

\item{\bf 1.}
Any group can be considered as an $\Omega$-group with empty set of
operations $\Omega$.
\item{\bf 2.}
A ring is an $\Omega$-group with commutative addition and
the set of operations
$\Omega$ consisting of one binary operation (multiplication),
satisfying the distributivity law.
\item{\bf 3.}
An algebra over a fixed field $F$ can be considered as an
$\Omega$-group with commutative addition and the set of
operations~$\Omega$ consisting of one binary operation
(multiplication) and $|F|$ unary operations (multiplication by
scalars) satisfying the well-known laws.

\enditem
Let $\V$ be a variety of $\Omega$-groups and let
$w(x_1,\dots,x_t)$ be an element of the free
(in the variety $\V$) algebra $F_{\V}(x_1,\dots,x_t)$.
For additive normal subgroups
$A_1,\dots,A_t$ of an $\Omega$-group $G\in\V$, we define
$w(A_1,\dots,A_t)$ as the normal subgroup generated by
all elements of the form
$w(a_1,\dots,a_t)$, where $a_i\in A_i$.

An element $w(x_1,\dots,x_t)\in F_{\V}(x_1,\dots,x_t)$ is called
\emph{multilinear} if
$$
w(A_1,\dots,A_i+A_i',\dots,A_t)=w(A_1,\dots,A_i,\dots,A_t)+
w(A_1,\dots,A_i',\dots,A_t)
$$
for all $i=1,\dots,t$ and
any normal subgroups $A_1,\dots,A_i,A_i',\dots,A_t$ of any
$\Omega$-group of the variety $\V$.

The outer commutators are multilinear elements of the absolutely free
group. The multilinear (in usual sense) expressions are multilinear
elements of the free algebra over a field.

Let $\C$ be a class of normal subgroups of an
$\Omega$-group $G$ such that
\item{1)}
$\C$ is closed with respect to the images under %!!
surjective endomorphisms of the $\Omega$-group $G$, finite sums
and finite intersections;
\item{2)}
any subfamily ${\N}\subseteq\C$ of the class $\C$ contains a
finite subfamily ${\F}\subseteq{\N}$ such that
$$
\sum_{N\in\N}N=\sum_{N\in\F}N.
$$

\enditem
In this case, we say that $\C$ is a \emph{class of large
normal subgroups}.
A function $\CODIM\:\C\to\R$ is called
a \emph{\(generalized\) codimension} if it has the following
properties:
\item{0)}
$\CODIM N_1\le\CODIM N_2$ if $N_1\supseteq N_2$;
\item{1)}
$\CODIM \phi(N) \le \CODIM N$
for each subgroup $N\in\C$ and each surjective endomorphism of the
$\Omega$-group~$G$;
\item{2)}
$\CODIM(N_1\cap N_2)\le \CODIM N_1+\CODIM N_2$
for all subgroups $N_1,N_2\in\C$;
\item{3)}
in any family $\N$ of subgroups from the class $\C$, there exist
$r\le\max\limits_{N\in\N}\CODIM N+1$ subgroups $N_1,\dots,N_r$
%p nekhorosho !!russkii
such that
$$
\sum_{N\in\N}N=\sum_{i=1}^r N_i.
$$

If $G$ is an algebra over a field and the class $\C$ consists of
all subspaces or all ideals (one-sided or two-sided) of finite
codimension, then the usual codimension can be regarded as a
generalized codimension.

If $G$ is a group and the class $\C$ consists of all normal
finite-index subgroups,
then, as a generalized codimension, we can take the binary
logarithm of the index. If the class $\C$ consists of all normal subgroups
whose index is a power of a fixed prime $p$, then, as a codimension
of a subgroup $N$, we can take the rank of the quotient $G/N$
(Property~3 holds by virtue of the Burnside basis theorem).

The class of all normal subgroups such that the corresponding
quotient groups satisfy a fixed smallness property~$\P$ can also
be regarded as a class of large subgroups. However, no codimension
is defined in this case.

%%%%%%%%%%%%%%%%%%%%%%%%%%%%%%%%%%%%%%%%%%%%%%%%%%%%%%%%%%%%%%%%%%%%%%%%%%
\s 3.
Main theorem

In the proof of the main theorem, we follow the argument in
[KlM09] generalizing it to multi-operator groups.

\Lemma 1.
Suppose that $w(x_1,\dots,x_t)$ is a multilinear element of the free
$\Omega$-group of some variety $\V$,
$m$ is a positive integer, $G\in \V$ is an $\Omega$-group, and
$\N$ is a finite family of its normal subgroups such that
$$
w(\underbrace{N,N,\dots,N}_{m\rm\;times},G,G,\dots,G)=0
\quad\hbox{for all $N\in\N$}.
$$
Then
$$
w(\underbrace{\^N,\^N,\dots,\^N}_{m-1\rm\;times},\^G,\^G,\dots,\^G)=0,
\quad\hbox{where } \^N=\bigcap_{N\in\N}\!\!N \hbox{ and }
\^G=\sum_{N\in\N}\!\!N.
$$

\Proof
$$
w(\underbrace{\^N,\^N,\dots,\^N}_{m-1\rm\;times},\^G,\^G,\dots,\^G)=
w(\underbrace{\^N,\^N,\dots,\^N}_{m-1\rm\;times},
\sum_{N\in\N}\!\!N,\^G,\dots,\^G)=
\sum_{N\in\N}\!\!w(\underbrace{\^N,\^N,\dots,\^N}_{m-1\rm\;times},
N,\^G,\dots,\^G).
$$
But $\^N\subseteq N$ and $\^G\subseteq G$;
therefore each term of the last sum is contained in the normal
subgroup
$$
w(\underbrace{N,N,\dots,N}_{m\rm\;times},G,G,\dots,G),
\quad\hbox{ which is trivial by assumption.}
$$

\noindent
As a corollary, we obtain the main theorem.

\proclaim{Main theorem}. Suppose that $G$ is an $\Omega$-group
belonging to a variety $\V$, $\C$ is a class of its large normal
subgroups, $w(x_1,\dots,x_t)\in F_\V(x_1,\dots,x_t)$ is a
multilinear element, $N\in\C$, and $w(N,\dots,N)=0$. Then $G$
contains an invariant under all surjective endomorphisms normal
subgroup $H\in\C$ satisfying the same identity $w(H,\dots,H)=0$.
In addition, if $\CODIM\:\C\to\R$ is a generalized codimension,
then
$$
\CODIM H\le f^{t-1}(\CODIM N),
$$
where $f^k(x)$ is the
$k$-th iteration of the function $f(x)=x(x+1)$.

\Proof
Let $\Ends G$ be the semigroup of all surjective endomorphisms of the
$\Omega$-group $G$.
Consider the normal subgroup
$
G_1=\!\sum\limits_{\varphi\in{\Ends G}}\!\!\!\!\!\varphi(N).
$
This subgroup is invariant under all surjective endomorphisms, is
large (i.e. belongs to $\C$), and $\CODIM G_1\le \CODIM N$ (if
$\CODIM$ is defined). Clearly, $G_1$ is the sum of a finite number
of the images of $N$ (because $N$ is large) and this finite number
does not exceed $\CODIM N+1$ (by the definition of the
codimension). Thus, $$ G_1=\sum_{k=0}^{p_1}{\varphi'_k(N)},
\quad\hbox{where $\varphi'_k \in \Ends G$ and $p_1\le l_0\:=\CODIM
N$}. $$ Now, consider the normal subgroup
$
N_1=\bigcap\limits_{k=0}^{p_1}{\varphi'_k(N)}.
$
Clearly, this subgroup is large too.
By Properties 1) and 2) of the codimension $\CODIM$, we have
$$
l_1\:=\CODIM N_1\le (p_1+1)\CODIM N=
(p_1+1)l_0\le (l_0+1)l_0=f(l_0).
$$
According to Lemma 1,
$$
w(N_1,\dots,N_1,G_1)=0.
$$

Similarly, we construct the large normal subgroups
$$
G_2=\!\sum_{\varphi\in{\Ends G}}\!\!\!\!\!{\varphi(N_1)}=
\sum_{k=0}^{p_2}{\varphi''_k(N_1)}
\quad\hbox{and}\quad
N_2=\bigcap_{k=0}^{p_2}{\varphi''_k(N_1)},
\quad\hbox{where }
\varphi''_k \in \Ends G
\hbox{ and }
p_2 \le \CODIM N_1= l_1 \le f(l_0).
$$
Clearly, $G_2$ is invariant under all surjective
endomorphisms of the $\Omega$-group $G$,
$$
\CODIM G_2 \le \CODIM N_1=l_1 \le f(l_0),
\hbox{ and }
l_2\:=\CODIM N_2\le (p_2+1)\CODIM N_1 =
(p_2+1)l_1\le f(l_1) \le f(f(l_0)).
$$
According to Lemma 1,
$$
w(N_2,\dots,N_2,G_2,G_2)=0.
$$

Continuing in the same manner, at the $t$-th step, we obtain in
$G$ an invariant under all surjective endomorphisms large normal
subgroup
$$
G_t=\!\sum_{\varphi\in{\Ends G}}\!\!\!\!\!{\varphi(N_{t-1})}=
\sum_{k=0}^{p_t}{\varphi_k^{(t)}(N_{t-1})},
\quad\hbox{where }
\varphi_k^{(t)}\in\Ends G.
$$
For this subgroup, we have
$$
w(G_t,\dots,G_t)=0 \quad\hbox{and}\quad \CODIM G_t \le
\CODIM N_{t-1}=l_{t-1} \le f(l_{t-2}) \le f(f(l_{t-3}))\le \dots\le
f^{t-1}(l_0).
$$
Thus, the subgroup $H=G_t$ is as required and the
theorem is proved.

\medskip
\noindent
Theorems 1 -- 4 are special cases of the main theorem:
$$
\matrix{
\hbox{Theorem 1:}&
\V=\{\hbox{Groups}\},\hfill&
\C=\{\hbox{Normal subgroups of finite index}\},\hfill&
\CODIM N =\log_2|G\!:\!N|;\hfill\cr\cr
\hbox{Theorem 2:}&
\V=\{\hbox{Algebras}\},\hfill&
\C=\left\{\vcenter{\hbox{\strut Subspaces or ideals}
\hbox{\strut of finite codimension}}\right\},
\hfill&
\CODIM=\codim;\hfill\cr\cr
\hbox{Theorem 3:}&
\V=\{\hbox{Groups}\},\hfill&
\C=\{N\nin G\;;\; G/N \hbox{ is a $p$-group}\},\hfill&
\CODIM N=\rank G/N;\hfill\cr\cr
\hbox{Theorem 4:}&
\V=\{\hbox{Groups}\},\hfill&
\C=\{N\nin G\;;\; G/N \hbox{ has the
property $\P$}\},\hfill&
\hbox{$\CODIM$ is not defined}.\hfill\cr\cr
}
$$

%%%%%%%%%%%%%%%%%%%%%%%%%%%%%%%%%%%%%%%%%%%%%%%%%%%%%%%%%%%%%%%%%%%%%%%%%%
\s 4.
Proof of Theorems $\bf 1'$ and $\bf 2'$

Theorem 1$'$ follows from the following proposition.%
\fn{%
To be more precise, Proposition 1 implies the assertion of Theorem
$1'$ on subgroups of bounded index. The finiteness of the total
number of finite-index maximal normal %!!russkii
subgroups with given identity follows from the proof of this
proposition. The relations between Proposition 2 and Theorem $2'$
are similar. }

\Proposition 1.
Suppose that a group $G$ contains sufficiently many normal subgroups
$N_1,\dots,N_m$ of index $n$ satisfying an
outer commutator identity $w(x_1,\dots,x_t)=1$
\(where $m$ is a sufficiently large number depending only on
$n$ and $t$\).
If
$$
\bigcap N_j \ne \bigcap_{j\ne k} N_j
\quad\hbox{for all $k=1,\dots,m$},
$$
then the group $G$ has a normal subgroup $X$
satisfying the same identity and
strictly containing one
of the subgroups~$N_j$ \(and, hence, the index of $X$ is strictly
less than $n$\).

\Proof
This proposition was, actually, proved in [KhM07a]. The subgroup
$X$ constructed in the proof of Proposition~1 of [KhM07a]
contains one of the subgroups $N_j$.

\medskip

Theorem $2'$ similarly follows from the following proposition
in~[KhM08].*$^{)}$

\Proposition 2 \rm([KhM08], Proposition 3).
{
Let $N_1,\dots,N_m$ be ideals of an algebra $A$ over a field $K$
and let $f(x_1,\dots,x_c)\in K\gp{x_1,\dots,x_c}$ be a multilinear
polynomial. Suppose that
\item{\rm (a)}
each ideal $N_i$ satisfies
the identity $f=0$ and $\dim A/N_i\le r$;
\item{\rm (b)}
$ \bigcap N_j \ne \bigcap\limits_{j\ne k} N_j
\quad\hbox{for all $k=1,\dots,m$}.$
\enditem
If $m\ge s(r,c)$ for some $(r,c)$-bounded number $s(r,c)$, then
there exists $k\in\{1,\dots,m\}$ such that the ideal
$N_k+\bigcap\limits_{j\ne k}N_j$ satisfies the identity $f=0$.
}

For the reader's convenience, we give an independent and simpler
proof of Theorem~$1'$.

%%%%%%%%%%%%%%%%%%%%%%%%%%%%%%%%%%%%%%%%%%%
\smallskip
\noindent
{\bf Proof of Theorem 1$'$.}
Let $\N$ be a set of finite-index subgroups of $G$ that
are maximal by inclusion among all normal subgroups satisfying the
outer commutator identity $w(x_1,\dots,x_t)=1$. We must prove that
$\N$ is finite.

If the family $\N$ is empty, then we have nothing to prove.
Otherwise,
consider a subgroup $G_0\in\N$. This subgroup satisfies the identity
$$
w_\sigma(G_0,\dots,G_0)=1
\quad\hbox{ for all } \sigma\in S_t.
\quad
%\hbox{ and for all } N \in\N.
\eqno (0) $$ Henceforth, $w_\sigma(x_1,\dots,x_t)$ denotes
$w(x_{\sigma(1)},\dots,x_{\sigma(t)})$, where $\sigma$ is a
permutation of degree $t$.

\smallskip

The subgroup $G_0$ has finite index. Therefore, the family of subgroups
$\{NG_0\;|\;N\in\N\}$ is finite and coincides with the family
$\{NG_0\;|\;N\in\N_1\}$, where $\N_1$ is a finite subfamily of the
family $\N$. The subgroup
$$
G_1=G_0\cap\bigcap\limits_{N\in\N_1}N
$$
has finite index and satisfies the equality
$$
w_\sigma(G_1,\dots,G_1,NG_0)=1
\quad\hbox{ for all } \sigma\in S_t
\quad
\hbox{ and for all } N \in\N.
\eqno (1)
$$
Indeed, by the choice of the family $\N_1$,
each product $NG_0$, where $N\in\N$, coincides with
a product $N_1G_0$ for some group
$N_1\in\N_1$ and $N_1\supseteq G_1\subseteq G_0$. Therefore,
$$
\eqalign{
w_\sigma(G_1,\dots,G_1,NG_0)=
w_\sigma(G_1,\dots,G_1,N_1G_0)=
&w_\sigma(G_1,\dots,G_1,N_1)w_\sigma(G_1,\dots,G_1,G_0)
\subseteq
\cr
\subseteq
&w_\sigma(N_1,\dots,N_1,N_1)w_\sigma(G_0,\dots,G_0,G_0)=1.
}
$$

\smallskip

The subgroup $G_1$ has finite index. Therefore, the family of subgroups
$\{NG_1\;|\;N\in\N\}$ is finite and coincides with the family
$\{NG_1\;|\;N\in\N_2\}$, where $\N_2$ is a finite subfamily
of the family $\N$. The subgroup
$$
G_2=G_1\cap\bigcap\limits_{N\in\N_2}N
$$
has finite index and satisfies the equality
$$
w_\sigma(G_2,\dots,G_2,NG_1,NG_1)=1
\quad\hbox{ for all } \sigma\in S_t
\quad
\hbox{ and for all } N \in\N.
\eqno (2)
$$
Indeed, by the choice of the family $\N_2$,
each product $NG_1$, where $N\in\N$, coincides with
a product $N_2G_1$ for some group
$N_2\in\N_2$ and $N_2\supseteq G_2\subseteq G_1\subseteq G_0$. Therefore,
$$
\eqalign{
&w_\sigma(G_2,\dots,G_2,NG_1,NG_1)=
w_\sigma(G_2,\dots,G_2,N_2G_1,N_2G_1)
=
\cr
=
&
w_\sigma(G_2,\dots,G_2,N_2,N_2)
w_\sigma(G_2,\dots,G_2,N_2,G_1)
w_\sigma(G_2,\dots,G_2,G_1,N_2)
w_\sigma(G_2,\dots,G_2,G_1,G_1)
\subseteq
\cr
\subseteq
&
w_\sigma(N_2,\dots,N_2,N_2,N_2)
w_\sigma(G_1,\dots,G_1,N_2,G_1)
w_\sigma(G_1,\dots,G_1,G_1,N_2)
w_\sigma(G_0,\dots,G_0,G_0,G_0).
}
$$
The first factor of the last product is trivial, because
the group $N_2$ satisfies the identity $w=1$. The second and
the third
factors are trivial by (1). The fourth factor
is trivial by (0).

\smallskip

Continuing in the same manner, we finally obtain a finite index
subgroup $G_{t-1}$ such that $$
w_\sigma(NG_{t-1},\dots,NG_{t-1})=1 \quad\hbox{ for all }
\sigma\in S_t \quad \hbox{ and for all } N \in\N. \eqno (t) $$ By
virtue of the maximality of all these subgroups $N$, this means
that $G_{t-1}\subseteq N$ for all $N\in\N$, i.e. $G_{t-1}\subseteq
\bigcap\limits_{N\in\N}N$ and, therefore, this intersection has
finite index. The finiteness of this index implies the finiteness
of the family $\N$, as required.

To obtain an estimate, it is sufficient to note that, if all subgroups
from the
family $\N$ have index not larger than~$n$, then
$$
|G{\,:\,}G_k|
\le
|G{\,:\,}G_{k-1}| n^{|\N_k|}
\quad\hbox{and}\quad
|\N_k|
\le
2^{|G{\,:\,}G_{k-1}|}
\hbox{ (this is a very rough estimate).}
$$
Therefore,
$$
|G{\,:\,}G_k|
\le
|G{\,:\,}G_{k-1}|\cdot n^{2^{|G{:}G_{k-1}|}},
\quad\hbox{i.e.}\quad
|G{\,:\,}G_{t-1}|
\le
F^{t-1}(n)
\quad\hbox{and}\quad
|\N|\le 2^{F^{t-1}(n)},
$$
where $F^k(x)$ is the $k$th iteration of the function $F(x)=xn^{2^x}$.

%%%%%%%%%%%%%%%%%%%%%%%%%%%%%%%%%%%%%%%%%%%%%%%%%%%%%%%%
\s 5.
Meta-virtually-solvable groups

In this section, we prove Theorem 5, i.e. we obtain the best
possible estimate on the `virtual derived length' of a group~$G$
containing a normal virtually solvable of derived length $s$
subgroup $A$ such that the quotient $G/A$ is virtually solvable of
derived length $t$. The second ``virtually" can easily be removed.
Indeed, replacing the group $G$ by its finite-index subgroup (the
preimage of the solvable finite-index subgroup of the quotient
$G/A$), we can assume  that the quotient $G/A$ is  solvable of
derived length $t$.

Next, by Theorem 1, the solvable of derived length $s$
finite-index subgroup $N$ of $A$ can be assumed to be
characteristic in $A$ and, hence, normal in $G$. To prove Theorem
5, it remains to show that the quotient $H=G/N$ contains a
finite-index subgroup which is solvable of derived length $\le
t+1$. But $H$ is an extension of the finite group $K=A/N$ by the
solvable of derived length $t$
%!!!russkii
group $G/A$.

A finite group $K$ has only finite number of automorphisms.
Therefore, the centralizer of this group in $H$ is of finite
index. Thus, passing to a finite-index subgroup, we can assume
that $K$ is contained in the centre of $H$. The central quotient
of $H$ has derived length $t$ and, therefore, $H$ itself is
solvable of derived length  $\le t+1$, as required.

\smallskip

It is easy to see that the estimate in %!!!
Theorem 5 cannot be improved. Indeed, consider, e.g., the central
product~$G$ (i.e. the direct product with amalgamated centres) of
infinitely many copies of the quaternion group of order 8. This
group $G$ is an extension of its finite (central) subgroup (of
order 2) by the elementary abelian 2-group of infinite rank. It is
easy to verify that this group has no abelian subgroups of finite
index. (To show this, one can use Theorem 1 once again: if there
is an abelian finite-index subgroup, then there exists a
characteristic abelian finite-index subgroup). %!!!??????NADO LI ETO??

Thus, an extension of a finite group (i.e. a virtually trivial
group, or a virtually solvable of derived length zero group) by an
abelian group is not necessarily virtually abelian. This example
can be modified to construct an extension of a virtually abelian
group by an abelian group which is not virtually metabelian.
Indeed, take
a faithful
%an exact
%??
complex representation $\phi\:G\to\GL(V)$
(e.g., the regular one) of the group $G$ described above. The
corresponding semidirect product $G_1=V\leftsemitimes G$ is an
extension of the virtually abelian group
$A=V\leftsemitimes\{\pm1\}$ by an (elementary) abelian group
$G_1/A\iso G/\{\pm1\}$. Suppose that $H$ is a finite-index
subgroup of $G_1$. Let us show that $H$ cannot be metabelian.
Indeed, $H$ must contain $V$ (because $V$, being a complex vector
space, has no proper finite-index subgroups). The quotient
$G_1/V\iso G$ contains a finite-index subgroup $H/V$. Therefore,
$H/V$ is nonabelian (since $G$ has no finite-index abelian
subgroups). Therefore, the commutator subgroup $\(H/V\)'$ contains
a nontrivial element~$g$ of order 2. Hence, the commutator
subgroup $H'$ of $H$ contains an element $x=ug$ (for some $u\in
V$) and all elements of the form $$ [x,v]=xvx^{-1}(-v)=\phi(g)v-v,
\quad\hbox{where $v\in V$}. $$ Since the representation $\phi$ is
faithful,
%??
%exact,
the space $V$ contains a vector $v$ not lying in the kernel
of the operator $\phi(g)-\id$. Thus, $H'$ contains a nonzero
vector $w=\phi(g)v-v$ such that $\phi(g)w=v-\phi(g)v=-w$. So,
$[x,w]=-2w\ne0$, i.e. $H'$ is nonabelian and $H$ is not
metabelian, as required.

Examples of higher derived lengths can be constructed in similar
fashion.

%%%%%%%%%%%%%%%%%%%%%%%%%%%%%%%%%%%%%%%%%%%%%%%%%%%%%%%%%%
\s 6.
The Burnside identity

In this section, we prove Theorem 6, i.e. we present a group virtually
satisfying the identity $x^p=1$, but having no large characteristic
subgroups of period $p$.  To construct such a group $G$, we use the
well-known technique dealing with periodic relations. We follow the book
[Olsh89]. A similar construction can be implemented on the base of the
book~[Adyan75].

Consider an infinite alphabet $X=\{a,x_1,x_2,\dots\}$ and the free group
$G(0)=F(X)$ with basis~$X$. The groups
$G(i)=\pres<X|R_i>$, where $i\ge1$, are defined inductively as follows.
Choose a set~$P_i$ of words of length $i$ in the group~$F(X)$ having the
following properties:
\item{1)}
no word in the set $P_i$ is conjugate in $G(i-1)$ to a power of a
shorter word;
\item{2)}
different words in the set $P_i$ are not conjugate in $G(i-1)$ to
each other and to the inverses of each other;
\item{3)}
the set $P_i$ is maximal (by inclusion) among all sets
satisfying conditions 1) and 2).

\enditem
The words in the set $P_i$ are called \emph{periods of rank $i$}.
Let us define a group
$
G(i)=\pres<X|R_i>
$
by setting
$$
R_0=\emptyset
\quad
\hbox{and}
\quad
R_i=R_{i-1}
\cup
\{u^{n_u}=1\;|\; u\in P_i\}
\hbox{ for $i\ge1$},
$$
where $n_u$ are some positive integers (depending on $u$).

Clearly, each nonidentity element of the group
$$
G=G(\infty)=\pres<X|\bigcup\limits_{i=1}^\infty R_i>
$$
is conjugated to a power of a period (of some
rank). It is known that, if all numbers $n_u$ are sufficiently large and
odd, then, in the group $G(\infty)$, the order of each period $u$ is
precisely $n_u$ (see, e.g., [Olsh89], Theorem 26.4).

Choose a sufficiently large prime $p$ and put
$$
n_u=\cases{
p& if $\phi(u)\notin\gp a \setminus\1$;\cr
p^2& if $\phi(u)\in\gp a \setminus\1$,\cr
}
$$
where $\phi\:F(X)\to\gp{a}_p\times \gp{x_1}_p\times\gp{x_2}_p\times\dots$
is the natural homomorphism of the free group onto the elementary abelian
$p$-group. Thus, in the group $G=G(\infty)$, the order of each period $u$
is either $p$ or $p^2$ depending on the value of $\phi(u)$.

\Lemma 2.
If the prime $p$ is sufficiently large, then
\item{\rm1)}
the group $G$ is periodic of exponent $p^2$;
\item{\rm2)}
the homomorphism $\phi$ induces a homomorphism
{\rm(denoted by the same letter)} of the group
$G$ onto the elementary abelian $p$-group;
\item{\rm3)}
an element $g$ of $G$ has order $p^2$ if and only if
$\phi(g)\in\gp a \setminus\1$ (the orders of other nonidentity elements
are $p$);
\item{\rm4)}
for each positive integer $i$ and each integer $k$, the
mapping $f_{i,k}\:X\to G$
that maps letter $x_i$ to $a^kx_i$ and fixes the other
letters of the alphabet $X$ extends to an automorphism of $G$.

\Proof
The first assertion is true, because
each element is conjugate to a power of a period and the periods have
orders $p$ and $p^2$. The second assertion follows from
the form $u^p=1$ of the defining relations.

To prove the third assertion, consider an element
$g$ of the group $G$. This element is conjugate to a power of a period:
$g=t^{-1}u^kt$. If $p$ divide $k$, then the equality $u^{p^2}=1$
(valid for any period~$u$) implies that the order of $g$ is either
$p$ or 1. On the other hand, $\phi(g)=\phi(u^k)=1$. So, in this case,
the assertion 3) is true.
If $p$ does not divide $k$, then
the order of $g$ coincides with the order of the period~$u$. On the other
hand, the inclusion $\phi(g)\in\gp a\setminus\1$ is equivalent to
the inclusion $\phi(u)\in\gp a\setminus\1$ and the assertion
follows from the above remark on the orders of periods.

Let us prove the fourth assertion. Note that it is sufficient to
show that the mappings under consideration extend to
endomorphisms, because, if this is true, then the endomorphisms
$f_{i,k}$ and $f_{i,-k}$ are mutually inverse.

To verify that the mapping $f_{i,k}$ can be extended to an
endomorphism, it is sufficient to show that $f_{i,k}$ transforms
the defining relations into valid equalities in $G$. Consider a
defining relation $u^l=1$, where $l$ is either $p$ or $p^2$
depending on the value of $\phi(u)$. The mapping $f_{i,k}$ induces
an automorphism of the elementary abelian $p$-group that fixes
each element of the subgroup $\gp a$. Therefore,
$\phi(f_{i,k}(u))\in\gp a\setminus\1$ if and only if
$\phi(u)\in\gp a\setminus\1$. Thus, by assertion 3), the element
$f_{i,k}(u)$ has the same order $l$ as the period $u$. Hence, the
mapping~$f_{i,k}$ transforms all the defining relations into valid
equalities and Lemma 2 is proved.

\medskip

\noindent{\bf Proof of Theorem 6.}
The group $G$ constructed above contains the subgroup
$N=\gp{x_1,x_2,\dots}\ker\phi$ of index~$p$. According to Lemma 2(3), this
subgroup satisfies the identity $x^p=1$.

Now, consider a characteristic subgroup $H$ of $G$. This subgroup
is invariant, in particular, under the automorphisms~$f_{i,k}$.
The automorphisms $f_{i,k}$ induce automorphisms of the elementary
abelian $p$-group
$$
E=\gp{a}_p\times\gp{x_1}_p\times\gp{x_2}_p\times\dots~.
$$
Therefore, the image
$\phi(H)$ of $H$ is a subgroup of $E$ invariant under all
automorphisms $f_{i,k}$. But any such invariant subgroup of~$E$ is
either trivial or containing $a$. Therefore, $H$ is either
contained in $\ker\phi$ (and, hence, the quotient~$G/H$ is
infinite and even has both-side infinite chains of normal
subgroups) or containing an element of order~$p^2$ (by Lemma~2(3))
and $H$ does not satisfy the identity $x^p=1$. This completes the
proof of Theorem 6.

%%%%%%%%%%%%%%%%%%%%%%%%%%%%%%%%%
\REFERENCES

\[Adyan75]
\thinspace
Adyan S.I.
The Burnside problem and identities in groups.
Moscow: ``Nauka", 1975.

\[BeK03]
Belyaev V.V., Kuzucuo\u glu M.
Locally finite barely transitive group
{//Algebra i Logika.} 2003. V.42.  no.3. P.261--270.

\[KaM82]
Kargapolov M.I., Merzlyakov Yu.I.
Fundamentals of group theory.
Moscow: ``Nauka", 1982.

\[KlM09]
Klyachko Ant. A., Melnikova Yu.B.
A short proof of the Khukhro--Makarenko theorem on large characteristic
subgroups with laws
{// Mat. Sbornik} (to appear).
See also
arXiv:0805.2747 .

\[Kur62]
Kurosh A.G.
Lectures on general algebra.
Moscow: ``Fiz.-Mat.Lit.", 1962.

\[Olsh89]
Olshanskii A.Yu.
Geometry of defining relations in groups.
Moscow: ``Nauka", 1989.

\[BrNa04]
Bruno B., Napolitani F.
A note on nilpotent-by-\v Cernikov groups
{// Glasgow Math. J.} 2004. 46, 211-215.

\[But08]
Button J.O.
Largeness of LERF and 1-relator groups
{// arXiv:0803.3805} , 2008.

\[KhM07a]
\thinspace
Khukhro E.I., Makarenko N.Yu.
Large characteristic subgroups satisfying multilinear commutator identities
{// J. London Math. Soc.} 2007. {V.75}. no.3. P.635--646.

\[KhM07b]\thinspace
Khukhro E.I., Makarenko N.Yu.
Characteristic nilpotent subgroups of bounded co-rank and
automorphically-invariant ideals of bounded codimension in Lie algebras
{// Quart. J. Math.} 2007. {V.58}. P.229--247.

\[KhM08]
Khukhro E.I., Makarenko N.Yu.
Automorphically-invariant ideals satisfying multilinear identities,
and group-theoretic applications
{// J. Algebra} 2008. {V.320}. no.4. P.1723--1740.

\end